\documentclass[11pt]{amsart}

\usepackage{color}
\usepackage{amsfonts, amsmath, amsfonts, amssymb}
\usepackage{graphicx}
\definecolor{RedClr}{rgb}{1,0,0}
\definecolor{BlueClr}{rgb}{0,0,1}
\definecolor{TextColor}{rgb}{0,0,0.5}
\definecolor{Violet}{rgb}{0.5,0,1}
\definecolor{Bordeaux}{rgb}{1,0.3,0.4}
\begin{document}

%%%%%%%%%%%%%%%%%%%%%%%%%%%%%%%%%
\newtheorem{thm}{Theorem}[section]
\newtheorem{cor}[thm]{Corollary}
\newtheorem{lem}[thm]{Lemma}
\newtheorem{prop}[thm]{Proposition}
\theoremstyle{definition}
\newtheorem{defn}[thm]{Definition}
\theoremstyle{remark}
\newtheorem{rem}[thm]{Remark}
\numberwithin{equation}{section} \theoremstyle{quest}
\newtheorem{quest}[]{Question}
\numberwithin{equation}{section} \theoremstyle{prob}
\newtheorem{prob}[]{Problem}
\numberwithin{equation}{section} \theoremstyle{answer}
\newtheorem{answer}[]{Answer}
\numberwithin{equation}{section}
\theoremstyle{fact}
\newtheorem{fact}[thm]{Fact}
\numberwithin{equation}{section}
\theoremstyle{facts}
\newtheorem{facts}[thm]{Facts}
\numberwithin{equation}{section}
\newtheorem{conj}[]{Conjecture}
\newtheorem{cntxmp}[thm]{Counterexample}
\numberwithin{equation}{section}
\newtheorem{exmp}[thm]{Example}
\numberwithin{equation}{section}
\newenvironment{prf}{\noindent{\bf Proof}}{\\ \hspace*{\fill}$\Box$ \par}
\newenvironment{skprf}{\noindent{\\ \bf Sketch of Proof}}{\\ \hspace*{\fill}$\Box$ \par}

\newcommand{\convGH}{\raisebox{-0.5em}{$\stackrel{\longrightarrow}{\scriptstyle GH}$}}

%%%%%%%%%%%%%%%%%%%%%%%%%%%%%%%%%%%%%%%%%%%%%%%%%%%%%%%%%%%%%%%%%%%%%%%%%%%%%

\title[]{Metric Curvatures and their Applications 2: Metric Ricci Curvature and Flow}
\author{Emil Saucan}

\address{%Max Planck Institute for Mathematics in the Sciences, Leipzig, Germany} %
Department of Applied Mathematics, ORT Braude College, Karmiel and Electrical Engineering Department, Technion, Haifa}

\email{semil@braude.ac.il, semil@ee.technion.ac.il}%

%\thanks{Research partially supported by Israel Science Foundation Grants 221/07 and 93/11.}%
%\subjclass{Primary: 53C44, 52C26, 68U05; Secondary: 65D18, 51K10, 57R40.}%
%\keywords{Combinatorial surface Ricci flow, metric curvature}%

\date{\today}

\begin{abstract}
In this second part of our overview of the different metric curvatures and their various applications, we concentrate on the Ricci curvature and flow for polyhedral surfaces and higher dimensional manifolds, and we largely review our previous studies on the subject, based upon Wald's curvature. 
In addition to our previous metric approaches to the discretization of Ricci curvature, we consider yet another one, based on the Haantjes curvature, interpreted as a geodesic curvature. 
We also try to understand the mathematical reasons behind the recent proliferation of discretizations of Ricci curvature. 
Furthermore, we propose another approach to the metrization of Ricci curvature, based on Forman's discretization, and in particular we propose on that  uses our graph version of Forman's Ricci curvature.
\end{abstract}

\maketitle

\tableofcontents

%%%%%%%%%%%%%%%%%%%%%%%%%%%%%%%%%%%%%%%%%%%%%

\section{Introduction: Why Ricci Curvature?}

%There are many approaches to Ricci curvature (new ones seem to appear, especially after Perelman's far-reaching(seminal, celebrated) work on -- and settling of -- the P\'{o}incare Conjecture. For surfaces, by far the best known and successful one is the Discrete Ricci one of Chow and Luo, based on circle packings. 

The natural question that would probably arise in the reader's mind is why, after surveying in the first part of this paper \cite{S15} a variety of essential notions of classical (smooth) curvature, we dedicate a separate paper solely to the metric discretizations of Ricci curvature? - After all, till relatively recently, this notion was relegated to the inner pages of Rimenannian curvature textbooks and to the mathematical study of General Relativity \cite{Besse}. True, a tremendous  amount of interest has been raised by Hamilton's \cite{Ha1},\cite{Ha2} and Perelman's \cite{Per1},\cite{Per2} work on the Ricci flow, and its application to the settling Thurston's Geometrization Conjecture (hence to the Poincar\'{e} conjecture). 
In particular, Chow and Luo's approach to Discrete Ricci curvature \cite{CL}, based on circle packings, has proved to be not only theoretically interesting and important, but also extremely useful and influential in various applications -- see e.g. \cite{Gu-Yau}, \cite{GGL}, \cite{ZZGLG}. Therefore, a drive to obtain a purely metric (in the sense of the first paper in this series \cite{S15}) Ricci flow is only natural. Still, being driven only by the quest for a metric analogue of the Combinatorial Ricci flow and its possible applications would have been, from the Metric Geometry viewpoint, somewhat mercenary. 

Indeed, another motivation for the geometer in understanding and generalizing Ricci curvature to discrete settings stems from a different direction, namely from the observation that Ricci curvature controls the volume growth of a manifold, a fact whose wide range of implications has been pointed out first by Gromov \cite{Gr-carte} (see also the deep results of Cheeger and Colding \cite{CC}, as well as Fukaya \cite{Fuk}). This insight conducted researchers towards a generalized notion of Ricci curvature bounds for Riemannian manifolds with density and, more generally, for metric measure spaces, namely the so called {\it curvature-dimension condition} ${\rm CD}(K,N)$, that was obtained by Lott-Villani \cite{LV} and Sturm \cite{St}. (We should note that these authors obtained a number of very strong results that confirm the correctness of these definition as a synthetic version of the notion of Ricci curvature bounded from above, such as locality, generalized versions of Meyer's and Bishop-Gromov's Theorems, generalized Poincar\'{e} inequalities, and stability under Gromov-Hausdorff limits, to number just a few; see \cite{Vi} for further details.) Naturally, the quest for a suitable version of this notion for a typical discrete case followed shortly \cite{BS}. 
A second strategy towards defining a Ricci curvature for weighted manifolds and graphs has its roots in the work of Bakry, Emery and Ledoux \cite{BE}, \cite{BL}. It is derived from the so called {\it Bochner-Weitzenb\"{o}ck formula} (see for instance \cite{J}), that relates curvature to the classical (Riemannian) Laplace operator. (We shall see that a quite different framework for the discretization of  Riicci curvature is also derived from the Bochner-Weitzenb\"{o}ck formula \cite{Fo}.) As with the Lott-Villani and Sturm approach, within this framework one is able to formulate a curvature dimension condition. Furthermore, in contrast with the Lott-Villani and Sturm definition, the one of Bakry and Emery easily lends itself to adaptation to graphs, as demonstrated by the works of Yau, Jost and their collaborators \cite{LY}, \cite{B---Y}, \cite{JL}. (see also \cite{Mu}). 
Finally, a third approach to Ricci curvature for weighted manifolds, due to  Morgan and his students \cite{Mo}, \cite{CHHWSX} originates -- yet again -- from another paper of Gromov, namely \cite{Gr} (the ideas of Bakry and Emery also playing a significant role in its development). Of the approaches to Ricci curvature of weighted structures it is the most direct and, perhaps, with the most evident geometric vein. Although powerful, this method employs simple formulae that represent straightforward generalizations of classical ones. Even so, only one attempt has been made so far (at least to best of or knowledge) to employ this approach in an applicative context -- see \cite{LLWZ}.

However, yet another, more geometric, view of Ricci curvature, proved itself to be most useful and intuitive  for the geometrization of such discrete spaces as graphs. This approach, due to Ollivier \cite{Ol1}, \cite{Ol2}, captures, in the discrete setting, the quantization of the growth of infinitesimal balls that represents one of the fundamental properties of Ricci curvature. Also, akin to the method of Lott-Villani and Sturm, it is based upon the notion of {\it optimal mass transportation}. It confirmed to be a powerful, yet intuitive, tool for the understanding of geometric (as well as analytic) properties of graphs and related structures (see, for instance, \cite{BJL}, \cite{JL}, \cite{LR}). 
In consequence, Ollivier's Ricci curvature represents an excellent apparatus for the intelligence of complex networks, in their various aspects, be they communication \cite{NLGGS},  \cite{WJB}, biological \cite{Allen1}, economical \cite{Allen2} or transportation \cite{Sh+} ones.

Moreover, besides these somewhat interrelated generalizations and discretizations of Ricci curvature, two approaches, both quite different from the ones above, as well as distinct from each other need to be noted, namely the ones due, to Stone \cite{St1}, \cite{St2},  and Forman \cite{Fo}. Since these works are strongly related to our own research that we shall detail below, we shall discuss each of them in more detail in the sequel.

It would seem, however, that we only managed to replace the original question (e.g. "Why should we devote a separate paper solely to discretizations of Ricci curvature?"), with an equally frustrating one, namely whether the multiplication of the various discretizations and generalizations of the classical notions is due, as one might perhaps suspect, to the excessive zeal of researchers who might have wished to contribute their own ideas on an attractive and ``hot'' topic, or s there an intrinsic reason why a ``definitive'' notion of discrete Ricci curvature proves to be so evasive?

It turns out that, there exists, at least in part, quite a deep scientific justification for this variety of discrete Ricci curvatures and that, in some sense, there is no possible ``ultimate'' discretization. More precisely, one is confronted with the following (unfortunately not known or acknowledged enough) result, due to Bernig \cite{Bernig2}, Proposition 8.1: For the large class of so called {\it singular spaces} (that includes the polyhedral manifolds that constitute the object of our interest in this paper) -- see e.g. \cite{BM} for the technical definition and further details , the exists 
no continuous, intrinsic tensor-valued distribution (measure) generalizing the Ricci tensor. In fact, no such generalization of the Riemann tensor is possible, either. It follows,  therefore, that the quest for ``perfect'' discretization of these tensorial classical curvatures to the context of polyhedral meshes, so dear to Graphics and Imaging community,  is a forlorn one, and to such ends, one has to do only with as good as possible approximations at most. More importantly, from a mathematical point of view, the result above proves, once more, that when passing from the classical (smooth) case to a discrete one, there is no choice but trying to capture one of the essential properties of the classical notion and use it as the basis definition in the discrete setting. Let us note here (with no little chagrin) that, while for geometers this is a clear and accepted paradigm, the Applied Mathematics (Imaging and related fields included) community is still reluctant to accept it.

\section{Metric Ricci Curvature and Flow for $PL$ Manifolds} \label{sec:MetricRicci1+2}
In this section we present some novel developments regarding Ricci curvature and flow for polyhedral manifolds, intrinsically based upon the Wald metric curvature %introduced 
discussed in detail in the first part of this study of metric curvatures in modern, computer related applications \cite{S15} %\ref{subsec:WaldCurv}. 
We begin by introducing a metric Ricci flow for surfaces (stemming from our paper \cite{S13}) and some of its applications and continue by proposing  (based on \cite{GS}) a metric Ricci curvature for $PL$ (polyhedral) manifolds, and test it, from the viewpoint of Synthetic Differential Geometry by proving a fitting variant (or rather variants --- see below) of the Bonnet-Myers Theorem.

\subsection{Metric Ricci Flow for Surfaces}
We begin with what represents not only the simplest case, but also the most important one, in view of the aforementioned works of Chow-Luo and Gu, namely that of the Ricci flow for surfaces.

\subsubsection{Smoothings and Metric Curvatures} \label{sec:Smoothings}
We begin by noticing that the Gromov-Hausdorff convergence results mentioned in \cite{S15} may appear to many a bit too general, at least when considering their applications to Imaging and related fields. One would like -- and not only in the applied context above -- wish to obtain similar results where `` second order'' curvature is also preserved (i.e. well approximated), not just ``first order geometry'', that is distances and where the approximating spaces are themselves smooth. This is, indeed, a highly natural question, to which one would like to give better, easier to use answers than the one provided by the observations brought (with anticipation) %immediately after Proposition \ref{prop:epsilon-net-conv} and mainly
in Remark 2.16 of \cite{S15}. One such answer is given by a results of Brehm and  K\"{u}hnel \cite{BK} as well as \cite{S13}, where a sketch of proof and its relation to Proposition \ref{prop:Waldapprox} below .

While less elementary than the one adopted by  Brehm and  K\"{u}hnel, our approach is, perhaps, more natural for geometers and topologists, as well as for Imaging applications, since it uses only standard analytic tools. In addition, it applies %extends
to a larger class of surfaces. Moreover, it captures %encapsulates
better the meaning of the notion of Hausdorff convergence and its interplay %(relationship)
with curvature: Instead of building, as in the Brehm-K\"{u}hnel approach, the smooth surfaces from a set of ``standard elements'' (cylinders, etc.), we consider instead {\it smoothings} $S_m^2$ (see, e.g. \cite{Mun}).
%(A similar approach of approximating discrete structures by smooth ones is adopted also in theoretical physics \cite{FL1}, \cite{FL2}, the paradigm therein being that the structure of space-time at the smallest scales is, in fact, discrete and that classical models are smooth approximations of these structures.)
%
%\marginpar{\tiny \bf Trebuie?!?}
%

\begin{figure}[htb] \label{fig:smoothing}
\begin{center}
\includegraphics[scale=0.26]{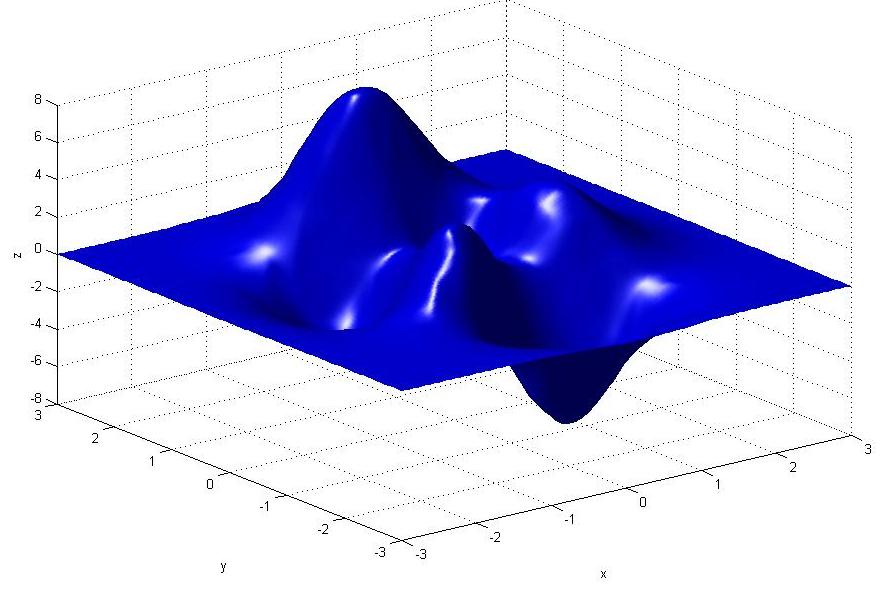}
\includegraphics[scale=0.26]{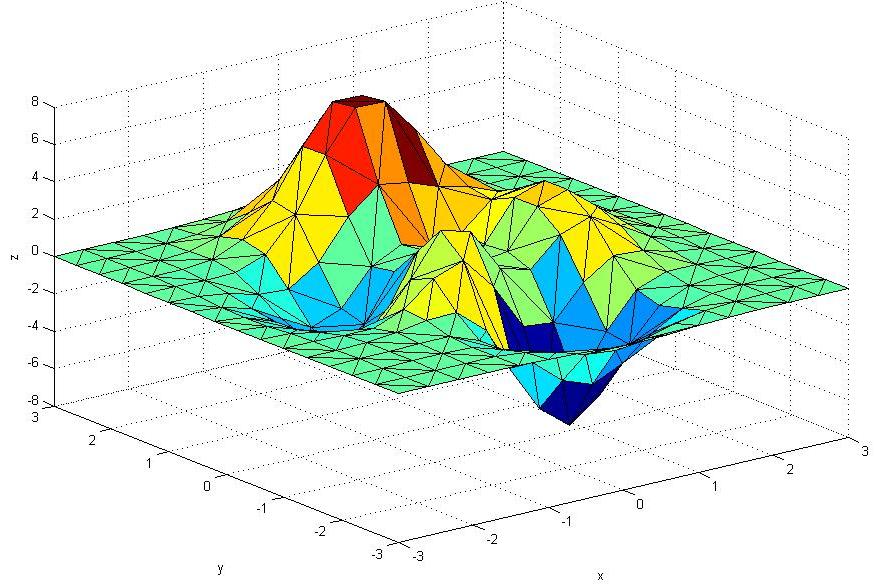}
%\pause
%\includegraphics[scale=.2]{Plouche.jpg}
%\pause
\end{center}
\caption{A smooth analytic surface (left) as the smoothing of one of its %own
$PL$ approximations (right).}
\end{figure}

Since, by \cite{Mun}, Theorem 4.8, such smoothings are $\delta$-{\it approximations}, that is metric approximations, and therefore, by \cite{Mun}, Lemma 8.7, they also are, for $\delta$ small enough, $\alpha$-{\it approximations}, i.e. angle approximations %\footnote{Recall that $\delta$-approximations are defined as follows: .....}
of the given piecewise-linear surface $S^2_{Pol}$, they approximate arbitrarily well both distances and angles on $S^2_{Pol}$. Therefore
%(it follows that)
angles, hence defects, are arbitrarily also well approximated. 
(For the precise technical definitions and proofs, see \cite{Mun}.)

The observations above amount, in fact, to a positive answer to the question --  not posed so far, to the best of our knowledge -- whether the metric curvature version of Brehm and H\"{u}hnels's basic result also holds, namely we have proved:

%\begin{quest}
%Given a combinatorial ($PL$) 2-manifold $S_{PL}$, does there exist a sequence of smooth surfaces (that is of class at least $C^2$) $S_m$, such that $K(S_m) \rightarrow K_W(S_{PL})$?
%\end{quest}

\begin{prop} [\cite{S13}]\label{prop:Waldapprox}
Let $S^2_{Pol}$ be a compact polyhedral surface without boundary. Then there exists a sequence $\{S^2_m\}_{m \in \mathbb{N}}$ of smooth surfaces, (homeomorphic to $S^2_{Pol}$), such that
\begin{enumerate}
\item
\begin{enumerate}
\item $S^2_m = S^2_{Pol}$ outside the $\frac{1}{m}$-neighbourhood of the 1-skeleton of $S^2_{Pol}$,
\item The sequence $\{S^2_m\}_{m \in \mathbb{N}}$ converges to % %\stackrel{\hspace{-0.5cm}\longrightarrow}{\mbox{$\scriptscriptstyle m \rightarrow \infty $}
$S^2_{Pol}$ in the Hausdorff metric;
\end{enumerate}
\item $K(S^2_{Pol}) \rightarrow K_W(S^2_{Pol})$, where the convergence is in the weak sense.
\end{enumerate}
\end{prop}

The result %proposition
above shows that, in fact, the metric approach to curvature is essentially equivalent to the combinatorial (angle-based) one, as far as polyhedral surfaces (in $\mathbb{R}^3$) are concerned.\footnote{This holds, of course, up to the specific type of convergence for the metric and combinatorial curvature, namely pointwise and in measure, respectively.}
In particular, as far as approximations of smooth surfaces in $\mathbb{R}^3$ are concerned, both approaches render,
in the limit, the classical Gauss curvature.) We should also stress that, in fact, the metric approach is more general, since it can be applied
%it applies
to a very large class of %(abstract)
metric spaces (see discussion at the end of this section).

\begin{rem}
The converse implication -- namely that Gaussian curvature $K(\Sigma)$ of a smooth surface $\Sigma$ may be approximated arbitrarily well by the Wald curvatures $K_W(\Sigma_{Pol,m})$ of a sequence of  approximating polyhedral surfaces $\Sigma_{Pol,m}$ -- is, as we have already mentioned above, quite classical.  %Indeed, .... -- see \cite{Bl}, ..?...
%(given the extensive treatment of the subject in the literature)
(For other approaches to curvatures convergence, see, amongst the extensive literature dedicated to the subject,  \cite{CMS} and \cite{BCM}, \cite{C-SM}, for the theoretical and applicative viewpoints, respectively.)
\end{rem}

\begin{rem}
The metric approach adopted here renders, in fact, a somewhat stronger result than that given in \cite{BK}, since no embedding in $\mathbb{R}^3$
%(or, for that matter, in any
is apriorily assumed, just in some $\mathbb{R}^N$ (even though as we have already noted above, this represents only a slight improvement).
%Moreover
Even more importantly, no change in the geometry of the 1-skeleton is made, not even in the neighbourhoods of the vertices.
\end{rem}

%----------------------------------------

\subsubsection{A Metric Ricci Flow}
From the  metric and curvature approximations results above, it follows that one can study the properties of the metric Ricci flow via those of its smooth counterpart, by passing to a smoothing of the polyhedral surface. The heavier machinery of metric curvature considered above pay off, in the sense that, by using it, the flow is purely metric and, moreover, the curvature at each stage (that is, for every ``$t$'') is given -- as in the classical context -- in an intrinsic manner, i.e. solely %purely
in terms of the metric.
\vspace*{0.25cm}

\paragraph{\it The Flow}
In direct analogy with the classical flow
\begin{equation} \label{eq:RicciSmooth}
\frac{dg_{ij}(t)}{dt} = - 2K(t)g_{ij}(t)\,.
\end{equation}
we define the {\it metric} Ricci flow as
\begin{equation} \label{eq:RicciComb}
\frac{dl_{ij}}{dt} = - 2K_il_{ij}\,,
\end{equation}
where $l_{ij} = l_{ij}(t)$ denote the edges (1-simplices) of the triangulation ($PL$ or piecewise flat surface) incident to the vertex $v_i = v_i(t)$, and $K_i = K_i(t)$ denotes the curvature at the same vertex. We shall discuss below in detail what proper notion of curvature should be chosen to render this type of metric flow meaningful.
%(In using proper, rather than partial derivatives notation, we follow \cite{CL}.)

%We are also interested(?!) in
We also consider %the close relative of (\ref{eq:RicciSmooth}),  %and its combinatorial version,
%namely
the normalized flow
\begin{equation} \label{eq:RicciFlow+}
\frac{dg_{ij}(t)}{dt} = (K - K(t))g_{ij}(t)\,,
\end{equation}
and its metric counterpart
\begin{equation} \label{eq:RicciComb+}
\frac{dl_{ij}}{dt} = (\bar{K} - K_i)l_{ij}\,,
\end{equation}
where $K,\bar{K}$ denote the average classical, respectively metric, sectional (Gauss) curvature of the initial %original
surface: $S_0$, $K = \int_{S_0}K(t)dA\big/\int_{S_0}dA$, and $\bar{K} = \frac{1}{|V|}\sum_{i=1}^{|V|}K_i$, respectively. (Here $|V|$ denotes, as usually, the number of the vertex set of $S_{Pol}$.)

Our approach to this type of metric flow should be evident in light of the previous discussion: Using smoothings, we approximate the metric Ricci curvature by its classical counterpart, thus replacing the metric flow with its smooth model.

\begin{rem}
Since for $PL$ surfaces (hence for their smoothings) $K_W(p) = 0$ for all points apart from vertices, the ensuing Ricci flow is stationary, except at vertices where the change rate is quite drastic. In this aspect, the metric Ricci flow introduced here resembles the combinatorial, rather than the smooth (classical) one. This should not be too surprising, given %the fact
that, as already stated, the initial motivation of considering the metric flow (and, a fortiori, smoothings) was to %smoothing is used just as a tool to obtain
gain a better understanding of some of the properties of the combinatorial flow.
\end{rem}

An advantage of the approach introduced here is that it automatically solves, by passing to infinitesimal distances, the asymmetry in equation \ref{eq:RicciComb}, that is caused by the fact that the curvature on two different vertices acts, so to say, on the same edge. We shall address this issue again shortly.

The main theoretical disadvantage of the approach embraced here is that it is not easy to extend it to higher dimensional manifolds. Indeed, it is not clear how to correctly define a flow for general high-dimensional manifolds (both for theoretical reasons and because of their applications, such as Medical Imaging,  Video, etc.), not least because 3-dimensional analogues of all the relevant results on the Ricci flow for smooth surfaces have yet to be obtained. Moreover, even defining a metric Ricci curvature in dimension 3 and higher is a non-trivial endeavor. See, however, the discussion and results in the following sub-section.

Therefore, developing  a purely metric Ricci flow -- that is one that does not make appeal to smoothings -- is highly desirable. One basic observation that has to be made is that one hast to take care of the lack of symmetry that we mentioned when we fist introduced the metric flow. From symmetry reasons, a natural way of defining the flow while addressing this is issue is (using the same notation was before):
\begin{equation} \label{eq:ImprovedRicciMetric}
\frac{d l_{ij}}{dt} = -\frac{K_i + K_j}{2} l_{ij}\;,
\end{equation}
where in this case, $K_i,K_j$ denote, of course, the Wald curvature at the vertices $v_i$ and $v_j$, respectively. (Note that, in fact, this expression appears also in the practical method of computing the combinatorial curvature, where it is derived via the use of a conformal factor \cite{Gu-Yau}.)

An important benefit of such a purely metric flow would be that, as in the case of combinatorial Ricci flow of Chow and Luo \cite{CL}, equation (\ref{eq:RicciSmooth})  becomes -- due to the fact that $K_i$ depends only the lengths of the edges $l_{ij}$, and not on their derivatives -- an  ODE, instead of a PDE, thence they are easier to study and enjoy better properties. In particular (and of importance in applications) the metric flow will have the backward existence property.

For further observations regarding the metric flow, as well as for some other possible improvements, see \cite{S13}.

\vspace*{0.25cm}

%----------------------------------

\paragraph{\it Existence and Uniqueness}
For the classical Ricci flow $dg_{ij}(t)/dt = 2K(t)g_{ij}(t))$ the (local) existence and uniqueness hold, on some {\em maximal} time interval $[0,T]; 0< T \leq \infty$ (see, e.g. \cite{CK}, as well as \cite{Ha}, for the original, different proof.\footnote{See also \cite{To}, Theorems 5.2.1 and 5.2.2 and the discussion following them for short exposition of the main steps of the proof.}) Moreover, the backward uniqueness of a solution (if existing) has been proven by Kotschwar \cite{Ko} (see also \cite{To} for a sketch of the proof).
Beyond the theoretical importance, the existence and uniqueness of the backward flow would allow us to find surfaces in the conformal class of a given circle packing (Euclidean or Hyperbolic).

More importantly, the use of purely metric approach, based on the Wald curvature (or any of other equivalent metric curvatures, for that matter), rather than the combinatorial (and metric) approach of \cite{CL}, allows us to give a first, tentative, purely theoretical at this point, answer to Question 2, p. 123, of \cite{CL}, namely whether there exists a Ricci flow defined on the space of all piecewise constant curvature metrics (obtained via the assignment of lengths to a given triangulation of 2-manifold).
Since, by Hamilton's results \cite{Ha} (and those of Chow \cite{Ch}, for the case of the sphere), the Ricci flow exists for all compact surfaces, it follows from the arguments above that the fitting metric flow exits for surfaces of piecewise constant curvature. In consequence, given a surface of piecewise constant curvature (usually a mesh with edge lengths satisfying the triangle inequality for each triangle), one can evolve it by the Ricci flow, either forward, as in the papers mentioned above, to obtain, after the suitable area normalization, the polyhedral surface of constant curvature conformally equivalent to it; or backwards -- if possible -- to find the ``primitive'' family of surfaces, including the ``original'' surface  conformally equivalent to the given one, (where, by ``original'', we mean the surface obtained via the backwards Ricci flow, at time $T$).

\begin{rem}
Note that is not necessarily true that all the surfaces obtained via the backwards flow are embedded (or, indeed, embeddable) in $\mathbb{R}^3$ --  a see the discussion below.
\end{rem}

\begin{prop}
Let $(S^2_{Pol},g_{Pol})$ be a compact polyhedral 2-manifold without boundary, having bounded metric curvature.
%\marginpar{\tiny Trebuie ``bounded''? - $M_{PL}$ compact!...}
Then there exists $T > 0$ and a smooth family of {\it polyhedral metrics}\footnote{see \cite{Hu}, \cite{S13}} $g(t), t \in [0,T]$, such that
\begin{equation}
\left\{
\begin{array}{ll}
\frac{\partial g}{\partial t} = -2K(t)g(t) & t \in [0,T]\,;\\
g(0) = g_{Pol}\,.
\end{array}
\right.
\end{equation}
(Here $K(t)$ denotes the Wald curvature induced by the metric $g(t)$.)

Moreover, both the forwards and the backwards (when existing) Ricci flows have the uniqueness of solutions property, that is, if $g_1(t), g_2(t)$ are two Ricci flows on $(S^2_{Pol}$, such that there exists $t_0 \in [0,T]$ such that $g_1(t_0) = g_2(t_0)$, then $g_1(t) = g_2(t)$, for all $t_0 \in [0,T]$.
\end{prop}

In fact, by combining our method with a result of Shi \cite{Sh}, we can extend the proposition above to complete polyhedral surfaces, as follows:

\begin{prop}
Let $(S^2_{Pol},g_{Pol})$ be a complete polyhedral surface, such that $0 < K_W \leq K_0$. Then there exists a (small) $T$ as above, such that there exists a unique solution of (\ref{eq:RicciSmooth}) for any $t \in [0,T]$.
\end{prop}

\vspace*{0.25cm}

\paragraph{\it Convergence Rate.}
A further type of result, highly important both from the theoretical viewpoint and for computer-driven applications, is that of the convergence rate. For the Ricci flow mostly applied in such s setting, namely for the combinatorial Ricci flow, it was proven in \cite{CL} that, in the case of background Euclidean (Theorem 1.1) or Hyperbolic (Theorem 1.2) metric, the solution -- if it exists -- converges, without singularities, exponentially fast to a metric of constant curvature. Using the classical results of \cite{Ha} and \cite{Ch}, since we already know that the solution exists and it is unique (see the subsection below for the nonformation of singularities), we are able to control the convergence rate of the curvature:

\begin{thm} \label{thm:flow-rate-conv}
Let $(S^2_{Pol},g_{Pol})$ be a compact polyhedral 2-manifold without boundary. Then the normalized metric Ricci flow converges to a surface of constant metric curvature. Moreover, the convergence rate is
\begin{enumerate}
\item exponential, if  $K < 0$; $\chi(S^2_{Pol}) < 0$;
\item uniform; if $K = 0$; %$\chi(S^2_{Pol}) = 0$; %
\item exponential, if $K > 0$. % $\chi(S^2_{Pol}) > 0$.
\end{enumerate}
\end{thm}

Recall that convergence rate of solutions is defined as follows:

\begin{defn}
A solution of (\ref{eq:RicciComb+}) is said to be {\it convergent} iff
\begin{enumerate}
\item $\lim_{t \rightarrow \infty}K_i(t) = K_i(\infty)$, for all $1 \leq u \leq |V|$, where $K_i(\infty) \in (0,2\pi)$;

\item $\lim_{t \rightarrow \infty}l_{ij}(t) = l_{ij}(\infty)$,  $l_{ij}(\infty) > 0$.
\end{enumerate}
A convergent solution is said to {\it converge exponentially fast} iff there exists constants $c_1,c_2$, such that, fora any $t \geq 0$, the following inequalities hold:
\begin{enumerate}
\item $|K_i(t) - \bar{K}_i| \leq c_1e^{-c_2t}$;

\item $|l_{ij}(t) - \bar{l}_{ij}| \leq c_1e^{-c_2t}$
\end{enumerate}
(The fitting definition for the flow (\ref{eq:RicciComb}) is immediate.)
\end{defn}

\begin{rem}
 A more realistic model for (gray-scale as well as color) images should be based on surfaces with boundary. Similar results can be obtained for this type of surfaces (see \cite{Br1}, \cite{Br2}), however we defer for further study %\cite{SAZ13a}
the detailed analysis, in this model, of the metric Ricci flow of images.
\end{rem}
\vspace*{0.25cm}

\paragraph{\it Singularities Formation.}
In the classical (smooth) case, by \cite{Ha}, Theorems 1.1 and 5.1, the Ricci flow evolves without singularities formation, even for surfaces of low genus.
Also, by \cite{CL}, Theorem 5.1, the combinatorial Ricci flow evolves without singularities on compact surfaces of genus $\geq 2$.
However, on surfaces of low genus, that are extremely important in Graphics and Imaging, singularities do form \cite{Gu}.

Combining our smoothing technique with Hamilton's classical results mentioned above, as well was with a more recent result of Topping \cite{To1}, we obtain the following results:

\begin{prop}
Let $(S^2_{Pol},g_{Pol})$ be a complete polyhedral 2-manifold, with at most a finite number of hyperbolic cusps (punctures), having bounded metric curvature and satisfying the noncollapsing condition below.

There exists $r_0 > 0$, such that, for all $x \in M$ the following holds:
\begin{equation} \label{eq:noncolapse}
{\rm Vol}_g\left(B_g(x,r_0)\right) \geq \varepsilon > 0.
\end{equation}
(Here, as usual, $B_g(x,r_0)$ denotes the open ball, in the metric $g$, of center $x$ and radius $r_0$.)

Then there exists a unique Ricci flow that contracts
%contracting
the cusps. Furthermore, the curvature remains bounded at all times during the flow.
\end{prop}
\vspace*{0.25cm}

%--------------------

\paragraph{\it Embeddability in $\mathbb{R}^3$.}
In Graphics, where one of the main problems solved  via the combinatorial Ricci flow is that of registration, by producing, via the flow, a conformal mapping from the given surface to one of the model surfaces  (see, e.g. \cite{Gu-Yau}, \cite{JKG}), the embeddabilty -- not necessarily isometric -- is both trivial and not of real interests. However, this aspect is highly significant in Image Processing (see \cite{ASZ11}), and, in fact, the results below were motivated precisely by these applicative aspects of the Ricci flow.

Here we mainly consider a problem regarding smooth surfaces, since, by now, the connection with the version for polyhedral surfaces is, we hope, quite clear. We should note that, by \cite{Mun}, Theorem 8.8, any $\delta$-approximation of an embedding is also an embedding, for small enough $\delta$. Since, as we have already mentioned, smoothing represent $\delta$-approximations, the possibility of using %passing
results regarding smooth surfaces to deduce facts regarding polyhedral embeddings is proven. (The reverse implication -- namely from smooth to $PL$ and polyhedral manifolds -- follows from the fact that the {\it secant approximation} is a $\delta$-approximation if the simplices of the $PL$ approximation satisfy a certain nondegeneracy condition -- see \cite{Mun}, Lemma 9.3.)

In the following $S_0^2$ denotes a smooth surface of positive Gauss curvature, and let $S_t^2$ denote the surface obtained at time $t$ from $S_0^2$ via the Ricci flow.

\begin{prop}
Let $S_0^2$ be the unit sphere $\mathbb{S}^2$, equipped with a smooth metric $g$, such that $K(g) > 0$. Then the surfaces $S^2_t$ are (uniquely, up to a congruence) isometrically embeddable in $\mathbb{R}^3$, for any $t \geq 0$.
\end{prop}

In fact, the result above can be slightly strengthened as

\begin{cor}
Let $S^2_0$ be a compact smooth surface. If $\chi(S_0^2) > 0$, then there exists some $t_0 \geq 0$, such that the surfaces $S^2_t$ are isometrically embeddable in $\mathbb{R}^3$, for any $t \geq t_0$.
\end{cor}

In contrast, for (complete) surfaces uniformized by the hyperbolic plane we have only a negative result:

\begin{prop}  \label{prop:non-embedd-hyp}
Let $(S^2_0,g_0)$ be a complete smooth surface, and consider the normalized Ricci flow on it. If $\chi(S^2) < 0$, then there exists some $t_0 \geq 0$, such that the surfaces $S^2_t$ are not isometrically embeddable in $\mathbb{R}^3$, for any $t \geq t_0$.
\end{prop}

For further related facts and comments, as well as some experimental confirmation of these results, see \cite{S13}.
\vspace*{0.25cm}

%-----------------------------------------

\subsection{Metric Ricci Curvature for $PL$ manifolds} \label{sec:MetricRicci2}

Our approach here stems from Formula (3.53) of \cite{S13} and it involves a combination of Wald's curvature and the earlier work of Stone \cite{St1}, \cite{St2}. The main (in fact essentially the only) part of Stones's work on which we rely upon is in his method of determining the relevant two sections and, of course, to decide what a direction at a vertex of a $PL$ manifold is. We begin our discussion of $PL$ Ricci curvature with this basic fact, as well as with related preliminaries. Next we explore various generalizations of the Bonnet-Myers Theorem, and we conclude by introducing a fitting version of scalar curvature.

\subsubsection{Definition and Basic Results}

In Stone's work, combinatorial Ricci curvature is defined both for the given simplicial complex $\mathcal{T}$, and also for its {\it dual complex} $\mathcal{T}^\ast$. In the later case, cells -- here playing the role of the planes in the classical setting of which sectional curvatures are to be averaged -- are considered. Unfortunately, his approach for the given complex, where one computes the Ricci curvature ${\rm Ric}(\sigma,\tau_1-\tau_2)$ of an $n$-simplex $\sigma$ in the direction of two adjacent $(n-1)$-faces, $\tau_1,\tau_2$, is less natural in a geometric context (even if useful in his purely combinatorial one), except for the 2-dimensional case, where it coincides with the notion of Ricci curvature in a direction (i.e., in this case, an edge). % -- see also Remark \ref{rem:defect} below).
On the other hand, passing to the dual complex will not confine us, since $(\mathcal{T}^\ast)^\ast = \mathcal{T}$ and, moreover -- and more importantly -- considering {\it thick} triangulations enables us to compute the more natural metric curvature for the dual complex and use the fact that the dual of a thick triangulation is thick, as we shall %present in
detail below. Moreover, working only with thick triangulations does not restrict us, however, at least in dimension $\leq 4$, since any triangulation admits a ``thickening'' -- see \cite{S05}.\footnote{This holds, as already mentioned, for any $PL$ manifold of dimension $\leq 4$, and in all dimensions for smoothable $PL$ manifolds, as well for any manifold of class $\geq \mathcal{C}^1$. Since the proof of one of our main results, regarding manifolds of dimension higher than 3, holds only for manifolds admitting smoothings, restricting ourselves only to such manifolds does not represent any additional obstruction.}

To define and compute the Ricci curvature of $\mathcal{T}$ and $\mathcal{T^*}$ and the connection between them, we have to %make appeal in an
essentially employ the thickness of the given complex.
Before proceeding further, it is imperative that we emphasize again the fundamental role of thickness in the sequel: Thickness ensures, by its definition, the fact that no degeneracy of the simplices occurs, hence no collapse and degeneracy of the metric can take place.
Moreover, %Reciprocally,
in its absence no uniform estimates for the edge lengths can be made, hence convergence of (dual) meshes and, as we shall see shortly, of their metric Ricci curvatures, can not be guaranteed.
Keeping in mind Munkres' definition of thickness \cite{S13} cf. \cite{Mun}, we begin by noting
% formula (\ref{eq:fat-Munkres}) --
that, since the length of the edge $l_{ij}^*$, dual to the edge $l_{ij}$ common to the faces $f_i, f_j$ equals $r_i + r_j$, the first barycentric subdivision\footnote{needed in the construction of the dual complex -- see e.g. \cite{Hu}}
of a thick triangulation is thick. (For planar triangulations, and also for higher dimensional complexes embedded in some $\mathbb{R}^N$, one can realize the dual complex (also in $\mathbb{R}^N$) by constructing the dual edges $l_{ij}^*$ orthogonal to the middle of the respective $l_{ij}$-s. To prove the thickness of the dual simplices, one has also to make appeal to the characterization of thickness in terms of dihedral angles (4.2) in \cite{S13}.
%(For the case of regular polytopes an explicit computation may be found in \cite{Be1}, Lemma 2.6.5.) {\tiny \bf -- Trebuie?!} }
%
To be sure, the notion of thickness  also makes sense for for general cells:

\begin{defn}
Let $\mathfrak{c} = \mathfrak{c}^k$ be a $k$-dimensional cell. The {\it thickness} (or {\it fatness}) of $\mathfrak{c}$ is defined as:
\begin{equation}
\varphi(\mathfrak{c}) = \min_{\mathfrak{b}}{\frac{{\rm Vol}(\mathfrak{b})}{{\rm diam}^l(\mathfrak{b})}}\,,
\end{equation}
where the minimum is taken over all the $l$-dimensional faces of
$\mathfrak{c}$, $0 \leq k$. (If ${\rm dim}\,\mathfrak{b} = 0$, then ${\rm
Vol}(\mathfrak{b}) = 1$, by convention.)
\end{defn}

The following facts now follow immediately:

\begin{lem}
Let $\mathcal{T}$ be a thick (simplicial) complex, and let $\mathcal{T^*}$ denote its dual. Then the following hold:

\begin{enumerate}
\item $\mathcal{T^*}$ is  thick.

\item Let $\delta(\mathcal{T}), \delta(\mathcal{T^*})$ denote the mesh of $\mathcal{T},\mathcal{T^*}$. Then
\begin{equation}
\lim_{\delta(\mathcal{T}) \rightarrow  0}(\mathcal{T}) = \lim_{\delta(\mathcal{T^*}) \rightarrow 0}(\mathcal{T^*})\,,
\end{equation}
where the convergence is in the Gromov-Hausdorff metric.
\end{enumerate}
\end{lem}

We can now return to the definition of Ricci curvature for simplicial complexes: Given a vertex $v_0$, in the dual of a $n$ dimensional simplicial complex, a {\it direction} at $v_0$ is just an oriented edge $e_1 = v_0v_1$. Since, %(by \cite{St1}?),
there exist precisely $n$ 2-cells, $\mathfrak{c}_1,\ldots,\mathfrak{c}_{n}$\,, having $e_1$ as an edge and, moreover, these cells form part of $n$ relevant variational (Jacobi) fields (see \cite{St1}), the Ricci curvature at the vertex $v$, in the direction $e_1$ is simply
\begin{equation} \label{eq:RicciCell}
{\rm Ric}(v) = \sum_{i=1}^nK(\mathfrak{c}_i)\,.
\end{equation}

Note that the index ``$i$'' in the definition (\ref{eq:RicciCell}) above runs from $1$, and not from $2$, as expected judging from the classical (smooth) setting.  This is due to the fact that we defined Ricci curvature by passing to the dual complex, with its simple but demanding (so to say) combinatorics. This fact has further implications -- see Theorem \ref{thm:ComparisonThm} and Remark \ref{rem:ComparisonThm} below.

To determine -- using solely metric considerations -- the sectional curvatures $K(\mathfrak{c}_i)$ of the cells $\mathfrak{c}_i$, we shall employ the ({\it modified}) {\it Wald curvature} $K_W$. Let us first note the role of the abstract open sets $U$ in the definition of Wald curvature (Definition 3.26 of \cite{S15}) is naturally played by the cells $\mathfrak{c}_i$. We can now formulate the following definition:

\begin{defn}
Let $\mathfrak{c}$ be a cell with vertex set $V_{\mathfrak{c}} = \{v_1,\ldots,v_p\}$. The ({\it metric}) {\it curvature} $K(\mathfrak{c})$ of $\mathfrak{c}$ is defined as:
\begin{equation}
K(\mathfrak{c}) = \min_{\{i,j,k,l\} \subseteq \{1,\ldots,p\}} %{1\leq i<j<k<l\leq p }
\kappa(v_i,v_j,v_k,v_l)\,.
\end{equation}
\end{defn}

\begin{rem}
In the definition above we presume that cells in the dual complex have at least 4 vertices. However, except for some totally degenerate (planar) cases, this condition always holds. However, even in this case Ricci curvature can be computed using a slightly different approach -- see the following remark.
\end{rem}

\begin{rem}
Note that by choosing to work with the dual complex we have restricted ourselves largely to considering solely submanifolds of $\mathbb{R}^N$, for some $N$ sufficiently large. However, in the case of 2-dimensional $PL$ manifolds this does nor represent  restriction, since, by a result of Burago and Zalgaller \cite{BZ} (see also \cite{S12b} and the references therein) 
%Huang, M. Keller, J. Masamune, R. Wojciechowski, A note on self-adjoint extensions of the Laplacian on weighted graphs})
 such manifolds admit isometric embeddings in $\mathbb{R}^3$, embeddings that, furthermore, are unique (up to isometries of the ambient space, of course).
\end{rem}

\begin{rem}
It should be emphasized that we have followed \cite{St1} only in determining the variational fields, but not in his definition of Ricci curvature. However, it is still possible (by dualization) to compute Ricci curvature according, more-or-less, to Stone's ideas, at least for the 2-dimensional case. (For more details see \cite{GS}.)
\end{rem}

%-----------------------

\subsubsection{A Bonnet-Myers Theorem}
Having introduced a metric Ricci curvature for $PL$ manifolds, one naturally wishes to verify that this represents, indeed, a proper notion of Ricci curvature, and not just an approximation of the classical notion. According to the synthetic approach to Differential Geometry (see, e.g.  \cite{Gr-carte}, \cite{Vi}), a proper notion of Ricci curvature should satisfy adapted versions of the main, essential theorems that hold for the classical notions. Amongst such theorems the first and foremost is Myers' Theorem (see, e.g., \cite{Be03}). And, indeed, fitting versions for combinatorial cell complexes and weighted cell complexes were proven, respectively, by Stone \cite{St1}, \cite{St2}, and Forman \cite{Fo}. Moreover, the Bonnet part of the Bonnet-Myers theorem, that is the one appertaining to the sectional curvature, was also proven for $PL$ manifolds, again by Stone -- see \cite{St3}, \cite{St0}.
\vspace*{0.25cm}

%----------------------------

\paragraph{\it The 2-dimensional case}
In the degenerate -- but of main importance in applications (see \cite{CL}, \cite{GY}, \cite{GS}) -- case of 2-dimensional manifolds, such a result is easy to formulate and prove, due to the fact that Ricci and sectional curvature essentially coincide:

\begin{thm}[Bonnet-Myers for $PL$ 2-manifolds -- Metric] \label{thm:BM-Metric}
Let $M^2_{PL}$ be a complete, connected 2-dimensional $PL$ manifold %without boundary,
such that

(i') There exists $d_0 > 0$, such that ${\rm mesh}(M^2_{PL}) \leq d_0$,
(where  ${\rm mesh}(M^2_{PL})$ denotes the mesh of the 1-skeleton of $M^2_{PL}$, i.e. the supremum of the edge lengths).

(ii') $K_W(M^2_{PL}) \geq K_0 > 0$.

Then $M^2_{PL}$ is compact and, moreover
\begin{equation} %\label{eq:BM}
{\rm diam}(M^2_{PL}) \leq \frac{\pi}{\sqrt{K_0}}\;.
\end{equation}
\end{thm}

\begin{rem}
Condition $(i)$, that ensures that the set of vertices of the $PL$ manifold is ``fairly dense''\footnote{in Stone's formulation (\cite{St0}, p. 1062).} is nothing but the necessary and common density condition for good approximation of both distances and of curvature measures, as we have already expounded in \cite{S13}, Section 4. The mere existence of such a $d_0$ is evident for a compact manifold, however it can't be apriorily be supposed for a general manifold, hence has do be postulated. In addition, to ensure a good approximation of curvature in secant approximation, this density factor has to be properly chosen (see, e.g. \cite{SAZ}), therefore tighter estimates for the mesh of the triangulation can be obtained using better curvature approximation. In this context we should also note that, apparently, the bound for diameter given by the proof above, is tighter than the one obtained by Stone in \cite{St3}, Theorem 3. Nevertheless, we should keep in mind that, in practice, one is more likely to encounter $PL$ surfaces as approximations of smooth ones.\footnote{and, obviously, $PL$ surfaces are $PL$ approximations of their own smoothings} However, the larger the mesh of the approximating surface (i.e. the ``rougher'' the approximation), the larger the deviation of the approximating triangles from the tangent planes (at the vertices), hence the more likely is to obtain large combinatorial curvature. Hence, there is a correlation between size of the simplices and curvature, even though not a straightforward one. No less importantly, an adequate choice of the vertices of the triangulation, also ensures, via the thickness property, the non-degeneracy of the manifold (and of its curvature measures), as we have detailed in Section 4 of \cite{S13}.
\end{rem}

\begin{rem} \label{rem:Kcomb}
Using the same approach a similar result for the combinatorial (defect) Gauss curvature of $M^2_{PL}$,
$%\begin{equation}
K_{Comb}(v_i) = 2\pi - \sum_{p=1}^{m_i}\alpha_p(v_i)\,
$%\end{equation}
 -- see \cite{GS}, where a different method of proof is also considered.
\end{rem}

\begin{rem}
The result above extends easily to polyhedral manifolds, since Wald curvature does not take into account the number of sides of the faces incident to a vertex, but only their lengths.
\end{rem}

\paragraph{\it The $n$-dimensional case, $n \geq 3$}
The proof of above does not extend immediately to higher dimensions, since, in general, no smoothing of a $PL$ manifold exists in dimension higher than $n > 4$ and, even if it exists, it is not necessarily unique, for $n \geq 4$ -- see \cite{Mun1}. However, if such a smoothing exists, the proof of does extend to any dimension, and we obtain the following $PL$ (metric) versions of the classical results:

    \begin{thm}[$PL$ Bonnet-Myers -- metric]%[Bonnet-Meyers for $PL$ manifolds(surfaces) -- Metric]
    \label{thm:BM-Metric+}
    Let $M^n_{PL}$ be a complete, $n$-dimensional $PL$, smoothable manifold without boundary, such that

(i'') There exists $d_0 > 0$, such that ${\rm mesh}(M^n_{PL}) \leq d_0$;

(ii'') $K_W(M^n_{PL}) \geq K_0 > 0$\,,

where $K_W(M^n_{PL})$ denotes the sectional curvature of the ``combinatorial sections'', i.e. the cells $\mathfrak{c}_i$.

Then $M^n_{PL}$ is compact and, moreover
\begin{equation} %\label{eq:BM}
{\rm diam}(M^2_{PL}) \leq \frac{\pi}{\sqrt{K_0}}\;.
\end{equation}
    \end{thm}

\begin{rem}
An approach similar to the one used in the proof above was also employed by Cheeger \cite{Ch} in a rather similar context. (It should be emphasized here that, as a byproduct of the results in \cite{GS} (and in this section), we also address -- using our own methods -- a problem posed by Cheeger in \cite{Ch}, Remark 3.5.)
\end{rem}

Before passing further on, we should underline the fact that, if we adopt the viewpoint of $PL$ (secant) approximations of smooth manifolds, then a number of problems arrise. Indeed, as we have already emphasized in Section 4 of \cite{S13}, even when such a smoothing $M^n$ ($n \geq 3$) exists, it is not probable that its sections provided  by $M^n_{PL}$, suffice to approximate well enough -- let alone reconstruct -- the Ricci curvature of $M^n$. Simply put, ``there are not enough directions'' in $M^n_{PL}$ to allow us to infer from the metric curvatures of a $PL$ approximation, those of a given smooth manifold $M^n$ (in fact, not not even a good approximation). On the other hand, increasing of the number of directions, i.e. of 2-dimensional sections (simplices) generates a decrease of the the precision of the approximation, due to the (possible) loss of thickness of the triangulation -- a problem which we have also discussed in detail in \cite4, Section 4.
The remarks above indicate that, unfortunately, in higher dimensional case, no general analogue of Myers' Theorem for $PL$ manifolds can be obtained by applying solely smoothing arguments). It is true that {\em a} Ricci curvature of the smooth manifold $M^n$ is obtained in terms of that of $M^n_{PL}$. However, it is not clear, in view of the paucity of sectional directions (i.e. possible 2-sections), how precisely is this connected to its discrete counterpart. Therefore, we can obtain, at best, an approximation result (with limits imposed by the thickness constraint -- see discussion above).

Due to this problematic aspect of the approximations approach, and in accord with the general mathematical drive towards generality, we adopt another approach in the following section.
\vspace*{0.25cm}

%-------------------

\paragraph{\it Alexandrov spaces}
We first pointing out the, perhaps not known well enough fact that Wald's curvature is essentially equivalent with the much more modern notion of {\it Alexandrov curvature}, at least for spaces in which there exists ``sufficiently many'' minimal geodesics (see, for instance, \cite{Pl}, Corollary 40), condition that certainly is fulfilled in $PL$ surfaces. Since Alexandrov curvature represents, by now, a quite classical and standard notion, and since introducing it formally here would take us too far afield, we will not bring here the technical definition and further details, but rather we refer the reader to, e.g. \cite{BBI}. However, we should mention that, in defining Alexandrov curvature, one makes appeal to {\it comparison triangles} in the model space (i.e. gauge surface $S_\kappa$), rather than quadrangles, as in the definition of Wald curvature. In fact, the sd-quads are, ``up to epsilon'' one of the (equivalent) ways of defining Alexandrov curvature. The reason for which we prefer to work with the Wald curvature, is that it is computable and, moreover, that it has even simpler, more practical approximations -- see \cite{SA09}).  (For the practical consequences of the similarities and differences between the two approaches, see \cite{S13}.)

It is, however, essential to notice that one has take into account the ``discrete'' nature of the types of spaces considered, hence to compute solely the Wald curvature of the 1-star neighbourhood of a vertex, as already stressed above, and not to consider (ever) smaller neighbourhoods, as perhaps natural in other contexts. This, however, agrees with the method of computing discrete curvature as angular defect, as employed, for instance, in \cite{BK} and in the Chow-Luo discrete Ricci flow \cite{CL} (as well as in many other instances -- see the bibliography for some of them). A positive consequence of his fact is that any such neighbourhood becomes a region having the same Alexandrov curvature bounded from below as the computed Wald one. Moreover, by the Alexandrov-Topogonov Theorem (see, e.g. \cite{Pl}, Theorem 43 and its proof, pp. 837-840), the whole surface becomes a space of curvature (Wald or Alexandrov) bounded from below.

In addition, taking into account only these ``discrete'' neighbourhoods is extremely important when equating the Wald and Alexandrov curvature, since it allows to avoid the blow-up of Alexandrov curvature at the vertices during smoothing. However, if one still wishes to consider smaller-and-smaller neighbourhood of the vertices (motivated, perhaps, by other applications then Imaging and Graphics, such as those in Regge calculus \cite{CMS}), one can resort to the basic approach of Brehm and K\"{u}hnel, that is ``rounding'' the edges by cylinders of radius $\varepsilon$ (without any change in curvature) and replacing the polyhedral cones at the vertices by smooth ``caps'', up to a predetermined  admissible error of, say, $\varepsilon_1$. Note that such a ``filtration'' of $K_W$ by Gaussian curvature (of the approximating smooth surfaces) is in concordance with common practices in Imaging, Vision and, indeed, in many applicative fields.
Moreover, considering only this ``discrete'' neighbourhoods is very important when equating the Wald and Alexandrov curvature, because it also allows us to avoid the blow-up of Alexandrov curvature at the vertices during smoothing.

By making appeal to the theory of Alexandrov spaces, and by using the equivalence of Wald and Alexandrov curvatures with the above mentioned provisos, a result of the desired type follows immediately from \cite{Pl}, Corollary 47, p. 840 and from the fact that $M^n_{PL}$ is locally compact:

\begin{thm}[Bonnet-Myers -- Alexandrov Spaces]
Let $M^n_{PL}$ be a complete, connected $PL$ manifold, such that $K_W(M^n_{PL}) \geq K_0 > 0$.

Then $M^n_{PL}$ is compact and, moreover
%
%\begin{equation} \label{eq:estimate1}
%{\rm diam}(M^2_{PL}) \leq D = D(K_0,d_0).
%\end{equation}
\begin{equation} %\label{eq:BM}
{\rm diam}(M^2_{PL}) \leq \frac{\pi}{\sqrt{K_0}}\;.
\end{equation}
\end{thm}

Unfortunately, determining weather a general $PL$ complex has Wald curvature bounded from below can be, in practice, quite difficult. However, in the special case of thick complexes (see definition in Section 1) one can determine a simple criterion as follows.

\begin{lem}
Let $M = M^n_{PL}$ be a complete, connected $PL$ manifold thickly embedded in some $\mathbb{R}^N$, such that $K_W(M^2) \geq K_0 > 0$, where $M^2$ denotes the 2-skeleton of $M$.
Then there exists $K_1 > 0$ such that $K_W(M^n_{PL}) \geq K_1 > 0$.
\end{lem}

%\begin{skprf}
%We indicate a proof the only for the case $n = 3$; the general case follows by a simple inductive argument.
%Consider an edge $e$ belonging to the 1-skeleton of $M^n_{PL}$. We have to show that $K_W(Q) \geq K_0 > 0$, for any quadruple incident to $e$.
%If $Q$ is one of the quadruples determined by the original cells of $M^2$ %(such as $Q_1$ in Figure 3)
%the condition is fulfilled trivially since $K_W(M^2) \geq K_0 > 0$.
%Otherwise the edges of $Q$ are either edges of the original cells, %(see Figure 3),
%or diagonals of such cells %(e.g. $d$ in Figure 3),
%or they connect vertices belonging to two different cells of the given complex. %(such as $\tilde{e}$ in Figure 3 connecting between vertices of the cells $\mathfrak{c}_2$ and $\mathfrak{c}_3$).
%But, it is quite standard to show that, by the fatness of the cells $\mathfrak{c}_i$, there exists a constant $c_1$ such that
%$\frac{1}{c_1}e \leq d \leq c_1e$. In a similar manner, using the boundedness from below of the angles as an equivalent definition of thickness, one can show that the fatness of the embedding implies that there exists a $c_2$ such that
%$\frac{1}{c_2}e \leq \tilde{e} \leq c_2e$.
%The desired conclusion follows from the two double inequalities above and from the continuity of the determinant function that defines the Wald curvature.
%\end{skprf}

The appropriate version of Bonnet-Myers now follows as a direct corollary:

\begin{thm}[Bonnet-Myers -- Thick Complexes]
Let $M = M^n_{PL}$ be a complete, connected $PL$ manifold thickly embedded in some $\mathbb{R}^N$, such that $K_W(M^2) \geq K_0 > 0$, where $M^2$ denotes the 2-skeleton of $M$.
Then $M^n_{PL}$ is compact and, moreover
%
%\begin{equation} \label{eq:estimate1}
%{\rm diam}(M^2_{PL}) \leq D = D(K_0,d_0).
%\end{equation}
\begin{equation} %\label{eq:BM}
{\rm diam}(M^2_{PL}) \leq \frac{\pi}{\sqrt{K_0}}\;.
\end{equation}
\end{thm}

\begin{rem}
We have formulated the theorem above in terms of piecewise flat manifolds since this is the case of most interest, both for theoretical ends, as well as application oriented ones. The most natural and useful instance in which such manifolds arise is that of secant approximations to smooth manifolds, as detailed and emphasized in the previous sections. However, the proof extends -- mutatis mutandis -- to the case of spaces whose simplices are modelled after spherical or hyperbolic spaces.
\end{rem}

\begin{rem}
This result is hardly surprising, given the fact that, by \cite{BGP}, Theorem 3.6, Myers' theorem holds for general Alexandrov spaces of curvature $\geq K_0 > 0$, and since, as already detailed above, Wald-Berestovskii curvature is essentially equivalent to the Rinow curvature, hence to the Alexandrov curvature. However,  instead of making appeal to the ``full force'' of the original proof, we gave, in the special case of $PL$ manifolds a simpler and more intuitive %garb, an alternative
proof of the Burago-Gromov-Perelman extension of Myers' Theorem.
\end{rem}

%-----------------------

\subsubsection{Metric scalar curvature}

Note that, while we defined and discussed into some depth Ricci curvature, we have not, up to this point, defined a fitting scalar curvature $K(\mathfrak{c})$ of a cell $\mathfrak{c}$.
Given our preceding discussion and keeping in mind the defining Formula (3.56) in \cite{S13}, the definition below is most natural:

\begin{defn}
Let $M = M^n_{PL}$ be an $n$-dimensional $PL$ manifold (without boundary). The {\it scalar metric curvature} ${\rm scal}_W$ of $M$ is defined as
\begin{equation} \label{eq:scal}
{\rm scal}_W(v) = \sum K_W(\mathfrak{c}),
\end{equation}
the sum being taken over all the cells of $M^*$ incident to the vertex $v$ of $M^*$.
\end{defn}

\begin{rem}
While the definition of scalar curvature of $M$ is defined, somewhat counterintuitively, by passing to its dual $M^*$, this approach is consistent with our approach to Ricci curvature (and also similar to Stone's original ideas -- see the discussion in 4.1 above).
\end{rem}

From this definition and our definition of sectional curvature of a cell, an immediate generalization of the classical curvature bounds comparison in Riemannian geometry follows (independently, in fact, of the chosen definition for the curvature of a cell):\footnote{Compare with \cite{Bernig1}, Theorem 1 and, for a measured version, with \cite{Bernig0}, Main Theorem 1.1.}

\begin{thm}[Comparison theorem] \label{thm:ComparisonThm}
Let $M = M^n_{PL}$ be an $n$-dimensional $PL$ manifold (without boundary), such that $K_W(M) \geq K_0 > 0$, i.e. $K(\mathfrak{c}) \geq K_0$, for any 2-cell of the dual manifold (cell complex) $M^*$. Then
\begin{equation} \label{eq:comp1}
K_W \lesseqqgtr K_0 \Rightarrow {\rm Ric}_W \lesseqqgtr nK_0\,.
\end{equation}

Moreover
\begin{equation} \label{eq:comp2}
K_W \lesseqqgtr K_0 \Rightarrow {\rm scal}_W \lesseqqgtr n(n+1)K_0\,.
\end{equation}
\end{thm}

\begin{rem} \label{rem:ComparisonThm}
\begin{enumerate}
\item
Inequality (\ref{eq:comp2}) can also be formulated in the seemingly weaker form:
\begin{equation} \label{eq:comp3}
{\rm Ric}_W \lesseqqgtr nK_0 \Rightarrow {\rm scal}_W \lesseqqgtr n(n+1)K_0\,,
\end{equation}
\hspace*{15cm}
\item
We should note that in all the inequalities above, the dimension $n$ appears, rather then $n-1$ as in the smooth, Riemannian case (hence, for instance one has in (\ref{eq:comp2}), $n(n+1)K_0$, instead of $n(n-1)K_0$\footnote{However, this holds even in dimension $n=3$!...} as in the classical case). This is a consequence of our definition (\ref{eq:RicciCell}) of Ricci (and scalar) curvature, that makes appeal to the dual complex of the given triangulation, therefore imposing standard and simple combinatorics, that allow only for such weaker bounds.\footnote{without affecting the analogue of the Bonnet-Myers Theorem -- see above.}
\end{enumerate}
\end{rem}

%--------------------------------------------------

%\subsection{Metric Ricci flow as Ricci flow on Alexandrov Surfaces}

\begin{rem}
Before concluding this section, we should mention the fact that we also discussed previously in \cite{S13}, namely that, by viewing compact $PL$ surfaces with as Alexandrov surfaces (with curvature bounded from below), one can apply the methods and results of Richard \cite{richard}, to obtain again the sought for results. Alternatively, one can mimic his approach, itself retracing the steps of the proof in the smooth case. While this is a task that we defer for further research, we take a deeper and more detailed look of the geometric consequences of viewing surfaces of bounded Wald curvature as Alexandrov surfaces in \cite{S14}.
\end{rem}

%-------------------------------------------------

\subsection{A Different Approach: From Geodesic Curvature to Ricci Curvature}

We briefly sketch below a different strategy of defining a notion of Ricci curvature for $PL$ and more general polyhedral  manifolds. To this end, we make appeal to the basic idea, stemming from Stone's work \cite{St1}, \cite{St2}, that we used in \cite{GS} and reviewed here in Section 2.2 above; in conjunction with Haantjes curvature in its role of geodesic curvature. The new, non-metric ingredient in our method is a discretization of the classical local Gauss-Bonnet theorem (see, e.g. \cite{doC}).

To begin with, less us notice that, in the case of piecewise flat manifolds (i.e the ``staple'' input for computations in Computer Graphics), the simplest and most direct approach is to view, for each triangle $T$ adjacent to an edge $e$, its {\it Menger curvature} $\kappa_M(T)$ (see \cite{S15}) as the sectional curvature of the section (plane) $T$. 
Since the Menger curvature represents a measure of the ``flatness'' of a triangle, $\kappa_M(T)$ gives a measure for the spreading of geodesics in the plane $T$, thus it provides, indeed, a discrete sectional curvature. 
Equipped with this metric version of sectional curvature, and using the same notations as in \cite{GS}, it is immediate to obtain the fitting versions of Ricci and scalar curvatures:

\begin{equation}
\kappa_M(e) \stackrel{\rm not}{=} {\rm Ric}_M(e) = \sum_{T_e \sim e}\kappa_M(T_e)\,,
\end{equation}  %\marginpar{\tiny \tt Recheck the notation}
where $T_e \sim e$ denote the triangles adjacent to the edge $e$; and

\begin{equation}
\kappa_M(v) \stackrel{\rm not}{=} {\rm scal}_M(v) = \sum_{e_k \sim v} {\rm Ric}_M(e_k) = \sum_{T \sim v}{\kappa_M(T)}\,;
\end{equation}
where,  $e_k \sim v, T \sim v$ stand for all the edges adjacent to the vertex $v$ and all the triangles $T$ having $v$ as a vertex, respectively. Note that ${\rm Ric}_M(e)$ captures, in accordance to the intuition behind $\kappa_M(T)$ the geodesic dispersion rate aspect of Ricci curvature. (See \cite{SSWJ} for a concise overview of the different aspects of Ricci curvature.)

While this manner of extending a simple metric curvature to the context of piecewise flat manifolds (and, in fact, of $PL$, in general), is quite direct and intuitive, it has also two limitations that hinder its usefulness in the study of such manifolds. 
The first such impediment is the fact that, being a discretization of the curvature notion for planar curves, Menger curvature of triangles is, intrinsically, always positive. Therefore, it can represent only a discretization of the absolute value of  sectional curvature (hence of Ricci and scalar curvatures as well). 
The second obstruction in extending this approach to more general types of meshes/manifolds resides in the fact that, Menger curvature, by its very definition, is restricted solely to triangles, thus is not applicable to more general cells.

However, making instead appeal to another type of metric curvature provides us with a solution for both problems. The second problem that we singled out above needs to be solved in the beginning, by considering instead of Menger curvature, {\it Haantjes curvature}. More precisely, given an edge $e = (u,v)$ and  a cell $\mathfrak{c}$, $\partial \mathfrak{c} = (u=v_0,v_1,\ldots,v_n=v)$, it is natural to define 

\[K_H(\mathfrak{c}) = 2\pi - \kappa_H(\pi)\,,\]
or rather, as we shall full justify below, as 

\begin{equation} \label{eq:K=-kH}
K_H(\mathfrak{c}) = 2\pi - \kappa_{H,e}(\pi)\,;
\end{equation}
where $\pi$ denotes the path $v_0,v_1,\ldots,v_n$, subtended by the chord $\bar{e} = \overline{v_0v_n}$, and where the notation $\kappa_{H,e}(\pi)$ emphasizes the fact that we compute the haantjes curvature of the path $pi$ subtended by the chord $e$. (Observe the slight change of notation conforming the one used in the definition of Haantjes curvature \cite{S15}.) With this definition we can also better extend the notion of, say, Ricci curvature, to $PL$ (piecewise flat) manifolds: By passing again (as in the previous section) to the dual manifold, one can ensure (except in the most degenerate case) that we deal with a polyhedral manifold, having 2-cells with more than 3 sides. This turns out to be essential here, even more than in Stone's approach, in defining a signed metric curvature, first of cells, thence of polyhedral manifolds. Note that we have to make appeal to Haantjes curvature, and not to the more simple, and well known, Menger curvature, due to the fact that, in general, there exists no natural, unique way of subdividing a cell (not embedded in some $R^n$) into triangles. 
%(See, however, the discussion at the end of this section.) 

The essential idea behind Formula (\ref{eq:K=-kH}) above is to make appeal to the local Gauss-Bonnet Theorem. Recall that, in the classical context of smooth surfaces it states that

\begin{equation}  \label{eq:SmoothGB}
\iint_DKdA + \sum_0^p\int_{v_i}^{v_{i+1}}k_gdl + \sum_{0}^p\varphi_i = 2\pi\chi(D)\,;
\end{equation}
where $D \simeq \mathbb{B}^2$ is a (simple) region in the surface $S^2$, having as boundary $\partial D$ a piecewise-smooth curve , of vertices (i.e. points where $\partial D$ is not smooth) $v_i, i = 1,\dots,p$, ($v_{n+1} = v_0$); $\varphi_i$ denotes the external angles of $\partial D$ at the vertex $v_i$; and $K$ and $k_g$ denote (as usually) the Gaussian and geodesic curvatures, respectively. 

We should note first that, in the absence of a background curvature, the very notion of angle is undefinable. Therefore, for abstract (non-embedded) cells, no meaningful (``decent'') notion of angle exists. In consequence, the last term on the left side of (\ref{eq:SmoothGB}) above has no signification, thus we should discard it. 
Indeed, the distances between non-adjacent vertices on the same cycle (apart from the path metric) are not defined, thus the third term in the left side of formula (\ref{eq:SmoothGB}) disappears. (We shall, however, reconsider it for embedded piecewise-flat manifolds.)

Moreover, for ``purely'' combinatorial $PL$ manifolds, the area of each cell can be prescribed (as it usually is) as being equal to 1. Moreover, one assumes (quite naturally) that curvature is constant on each cell. Therefore, the first term in the left side of (\ref{eq:SmoothGB}) reduces simply to $K$. % In the same context, there is no natural (or even ``decent'' notion of angle). 
In addition, given that $D$ is a 2-cell, thus $\chi(D) = 1$. It follows, that in such a setting we obtain
\begin{equation*}
K = 2\pi - \int_{\partial D}k_gdl\,.
\end{equation*}

It is tempting to next consider $\partial D$ as being composed of segments (on which $k_g$ vanishes), except at the vertices, thus rendering expression above as 
\begin{equation}
K_{H,e}(\mathfrak{c}) = 2\pi - \sum_1^{n-1}\kappa_{H,e}(v_i)\,;
\end{equation}

However, in abstract piecewise-flat manifolds, one can not define a non-trivial Haantjes curvature for each of the vertices since, as already noted above, no distance (apart from the one given by the path metric) between the vertices $v_{i-1}$ and $v_{i+1}$ can be considered. In fact, in this general case,  neither can the arc (path) $\pi = v_0v_1\ldots v_n$ be truly viewed as smooth. Therefore, in this context, we have no choice but to replace the right term in Formula (2.24) above by $\kappa_{H,e}(\pi)$, where it should be remembered that $\pi$ represents the path $v_0,v_1,\ldots,v_n$,  of chord $\bar{e} = \overline{v_0v_n}$. We have thus obtained the proposed Formula (\ref{eq:K=-kH}). Thus, the suggested definition for the Ricci curvature of a general metric cell-complex becomes

\begin{equation}
{\rm Ric}(e) = \sum_{\mathfrak{c} \sim e}\kappa_{H,e}(\mathfrak{c})\,,
\end{equation}
where the sum is taken over all the 2-cells $\mathfrak{c}$ adjacent to $e$.

Before concluding this section, we should retrace our steps in deriving the proposed formula for computing the discrete Ricci curvature based on a fitting adaptation of the local Gauss-Bonnet formula, and point out that, for piecewise-flat manifolds embedded in $\mathbb{R}^3$ (but not only -- see for instance \cite{Alsing+}, \cite{Alsing++}) one can establish a discrete version of Gaussian/sectional curvature much closer to the classical one.\footnote{Indeed, piecewise-flat surfaces (and $PL$ manifolds in general), can be viewed as bridging the gap between smooth manifolds and discrete (metric) spaces.} Since, in this case, proper angles at the vertices $v_i$ are defined, one can also compute the curvatures $\kappa_H(v_i)$ and, moreover, the area of each cell can be calculated (and it is independent on the specific subdivision of the cell into triangles). 
This new method of defining, via elementary computations, sectional curvature of cells for embedded polygonal manifolds is distinct both from Stone's approaches \cite{St1,St2} and \cite{St0}, and our previous one \cite{GS}. 
We defer the implementation, to the context of Graphics, of such work for later study. Moreover, by adapting the arguments of \cite{GS}, one can seemingly obtain similar results to those therein, namely regarding the convergence of this discretization of Ricci curvature, and analogues of the Bonnet-Myers and comparison theorems.

%However, piecewise flat manifolds, and in particular triangular and tetrahedral meshes in $\mathbb{R}^3$,

%=====================================

\section{An Alternative Approach: Forman's Ricci Curvature}

If one is willing to sacrifice the Wald metric curvature paradigm, an alternative approach -- albeit somewhat more abstract -- suggest itself, namely that based on Forman's Ricci curvature \cite{Fo}. While admittedly less geometric in nature (its development being based on constructing a Bochner Laplacian on weighted $CW$ complexes), it incorporates naturally the given metric (as, indeed, the idea residing behind its introduction is to discretize the metric and curvature of a Riemannian manifold).

Before proceeding further to the manner one can employ this curvature to the metric setting at hand, let us first give its definition:

Given a $p$-dimensional cell (or $p$-cell, for short) $\alpha = \alpha^p$, one can define the {\it curvature function}:

\begin{equation} \label{saucan-eqn:4}
\mathcal{F}_p = \langle F_p(\alpha),\alpha \rangle,
\end{equation}
where $F_p: C_p \rightarrow C_p$ is being regarded as a linear function on $p$-chains of cells (see \cite{Fo} for details), the scalar product appearing in the formula being defined (rather standardly) by $<\alpha,\alpha> = w, <\alpha,\beta> = , \alpha \neq \beta$,  and where 

\begin{equation} \label{eq:Forman-General}
\mathcal{F}(\alpha^p) = w(\alpha^p)\Big[\Big(\sum_{\beta^{p+1}>\alpha^p}\frac{w(\alpha^p)}{w(\beta^{p+1})}\;
%\]
%%
%\[
+ \sum_{\gamma^{p-1}<e_2}\frac{w(\gamma^{p-1})}{w(\alpha^p)}\Big)\; - %\]
\end{equation}
\[
\hspace*{0.5cm}
-\sum_{\alpha_1^p\parallel \alpha^p, \alpha_1^p \neq \alpha^p}\Big|\sum_{\beta^{p+1}>\alpha_1^p,\beta^{p+1}>\alpha^p}\frac{\sqrt{w(\alpha^p)w(\alpha_1^p)}}{w(\beta^{p+1})}\: 
- \sum_{\gamma^{p-1}<\alpha_1^p,\gamma^{p-1}<\alpha^p}\frac{w(\gamma^{p-1})}{\sqrt{w(\alpha^p)w(\alpha_1^p)}}\Big|\:\;\Big]
\]
the notation $\alpha < \beta$ meaning that $\alpha$ is a face of $\beta$, and
the notation $\alpha_1 \parallel \alpha_2$ signifies that the
simplices $\alpha_1$ and $\alpha_2$ are {\it parallel}, parallelism being defined as follows:\\\\
{\bf Definition 3}  \label{saucan-def:parallel}
Let $\alpha_1 = \alpha_1^p$ and $\alpha_2 = \alpha_2^p$ be two
p-cells. $\alpha_1$ and $\alpha_2$ are said to be  {\em parallel}
($\alpha_1
\parallel \alpha_2$) iff either:
(i) there exists $\beta = \beta^{p+1}$, such that $\alpha_1, \alpha_2 <
\beta$; or (ii) there exists $\gamma = \beta^{p-1}$, such that
$\alpha_1, \alpha_2 > \gamma$ holds, but not both conditions simultaneously. 

In the special case $p = 1$, one is conducted, by analogy to the smooth case, to the following 

\begin{defn} \label{saucan-def:Ricci}
Let $\alpha = \alpha^1$ be a 1-cell (i.e. an edge). Then the {\em Forman Ricci curvature} of $\alpha$ is defined as:
\begin{equation}
{\rm Ric}_F(\alpha)= \mathcal{F}_1(\alpha).
\end{equation}
\end{defn}

Therefore, it follows from Formula (\ref{eq:Forman-General}) above that

\begin{equation} \label{eq:Forman-2d}
{\rm Ric}_{\rm F} (e) = \omega (e) \left[ \left( \sum_{e \sim f} \frac{\omega(e)}{\omega (f)}+\sum_{v \sim e} \frac{\omega (v)}{\omega (e)}	\right) \right. %\nonumber \\
- \left. \sum_{\hat{e} \parallel e} \left| \sum_{\hat{e},e \sim f} \frac{\sqrt{\omega (e) \cdot \omega (\hat{e})}}{\omega (f)} - \sum_{v 	\sim e, v \sim \hat{e}} \frac{\omega (v)}{\sqrt{\omega(e) \cdot \omega(\hat{e})}} \right| \right] \; .
\end{equation}
%
%(See also Figure... .)

The comprehensiveness of Formula \ref{eq:Forman-General} above naturally raises questions about its range of applicability and effectiveness in two regards: The possible choices of weights and the specific type of cell complex considered.

As far as the weights are considered, Forman's work assures us (Theorems 2.6 and 3.9) that any set of weights encountered in practice can be used (or at least arbitrarily well approximated). However, both in theory and in practice, {\it natural} (or {\it geometric}) sets of weights, i.e. proportional the dimension of the cell (s.a. length, area, volume) are preferable. This is not due just to the specific techniques employed in \cite{Fo}, it is in fact ingrained in the very idea resting behind Forman's article, namely that of studying curvatures (and Laplacians) of ``good'' discretizations, e.g. triangulations, cubulations, etc., of smooth Riemannian manifolds, hence of their metrics and volume elements. We have employed this  type of weights in the case of complexes consisting from square grids naturally arising in the Image Processing setting, in \cite{SAWZ}, \cite{ASZ11}, \cite{SSAZ}. 
However, much more general (positive) weights can be considered, as already noted above. This is extremely important in many practical applications, for instance in Medical Imaging, where, for instance, MRI images are given, basically, by proton densities,\footnote{Curvature measures also are natural occurring weights in Imaging and Graphics, and weights can be attached to textures as well.} but also in other fields, e.g. Manifold Learning where weights can be arise as probabilities or importance attached to a specific region, etc. We shall shortly return to the problem of ``correctly'' prescribing Ricci curvature to such weighted cell complexes.

Regarding the type of underlying complex, a number of observations are mandatory: First, there are two specific types of structures that warrant special attention, since they appear naturally and often in many instances, both practical and theoretical. 
The first such instance is that of {\it cube manifolds}, for which, due to the simple form of parallelism holding in their case, Formula \ref{eq:Forman-General} is drastically simplified -- see \cite{Fo}, Theorem 8.1 for the formula in the $n$-dimensional case, and the papers mentioned above for the weighted 2 (and 3) dimensional complexes encountered in Imaging.\footnote{See also \cite{S15a} for an application of metric curvatures to the problem of smoothability of metrics on cube complexes and the references therein for the theoretical setting and uses of such complexes.} 
For the weighted square and cube grids that arise in classical and volumetric Imaging these formulas become
\begin{equation} \label{eq:Ricci-Forman2D}
\hspace*{-0.5cm}{\rm Ric}_F(e_0) =
w(e_0)\left[\left(\frac{w(e_0)}{w(c_1)} +
\frac{w(e_0)}{w(c_2)}\right) -
\left(\frac{\sqrt{w(e_0)w(e_1)}}{w(c_1)} +
\frac{\sqrt{w(e_0)w(e_2)}}{w(c_2)}\right)\right]\,,
\end{equation}
and
\begin{equation}  \label{eq:Ricci-Forman3D}
{\rm Ric}_F(e_0) =
w(e_0)\left[w(e_0)\left(\sum_{1}^{4}\frac{1}{w(c_i)}\right) -
\sqrt{w(e_0)}\left(\sum_{1}^{4}\frac{\sqrt{w(e_i)}}{w(c_i)}\right)\right],
\end{equation}
respectively. (Here and above we took into account that, for digital images, the vertices' weights are always 0.)

The second important case is that of $PL$ manifolds, occurring (and, indeed, motivated) principally as triangulations and piece-wise flat approximations of smooth surfaces and higher dimensional manifolds, both in a theoretical setting, but most commonly in Graphics, CAD, Imaging and related fields. Here again a simplified Formula \ref{eq:Forman-General} can be written (due, in this case, to the lack of parallelism for the edges of the 2-cells adjacent to an edge $e$):

\begin{equation}
\label{eq:Forman-2d-tr}
{\rm Ric}_F(e) = \omega(e)\left[\sum_{t > e}\frac{\omega(e)}{\omega(t)} + \left(\frac{\omega(v_1)}{\omega(e)} + \frac{\omega(v_2)}{\omega(e)}\right) - \sum_{\tilde{e} \sim e}\left|\sum_{t > e,\, t > \tilde{e}}\frac{\sqrt{\omega(e)\omega(\tilde{e})}}{\omega(t)} - \sum_{v < e,\, v < \tilde{e}}\frac{\omega(v)}{\sqrt{\omega(e)\omega(\tilde{e})}}\right|\right]\,.
\end{equation}
Moreover, given that $e$ has only two 0-dimensional faces (i.e. the vertices $v_1,v_2$ we can replace the term

\[\sum_{t > e,\, t > \tilde{e}}\frac{\sqrt{\omega(e)\omega(\tilde{e})}}{\omega(t)}\]

with

\[\sum_{t > e,\, \tilde{e} > v_1}\frac{\sqrt{\omega(e)\omega(\tilde{e})}}{\omega(t)} + \sum_{t > e,\, \tilde{e} > v_1}\frac{\sqrt{\omega(e)\omega(\tilde{e})}}{\omega(t)}\,.\]

We should also note that, for combinatorial weights, the addition to the complex $X$ of each triangle $\bf t$ adjacent to an edge $e$ increases the Forman Ricci curvature by 3, %\marginpar{\tiny \bf Verifica daca intr-adevar 3!}
that is
\begin{equation}
{\rm Ric_F} (e | X+t) = {\rm Ric_F} (e | X) + 3 \; .
\end{equation}
(For details, see \cite{WSJ2}.)

However important these particular cases might be, it is still important to emphasize here again that, since Ricci curvature is a measure attached to edges, it is determined solely by the 2-skeleton, given that edges represent the common boundaries of the 2-cells of complex. This is clearly evident in the general formula for the Forman-Ricci curvature of an edge (\ref{eq:Forman-2d}).

Thus, Forman-Ricci curvature of even more complicated structures with 2-skeleta of the types above can be easily computed, independently of the structure of the higher dimensional cells. Experimental results so far indicate that, at least for images and triangular meshes, Forman's Ricci curvature approximates well, as a measure, the Gaussian curvature  with whom Ricci curvature essentially identifies in dimension $n = 2$. A formal analysis of the convergence properties of $\mathcal{F}(\alpha^p)$ (and of the associated Bochner Laplacian) are currently in progress \cite{LS}.

\begin{rem}
	As we have emphasized in \cite{SJ}, the Forman Ricci curvature has a further inherent advantage, namely that it allows for the understanding of the algebraic topological properties (homology and first homotopy groups), rendering it as a useful tool in the intelligence of discrete structures, instead -- or in complement of -- Persistent Homology.)
\end{rem}

In this context one could use the the vicinity (dual graph) of a given tessellation (e.g triangulation or cubulation), with node  weights that concentrate the weights of the respective dual faces, and compute the graph Forman-Ricci curvature and associated flow  developed for Complex Networks in \cite{SMJSS}, \cite{WJS}, \cite{WSJ1}. 
Recall that Forman's Ricci curvature for an edge is given by the following formula:

\begin{equation} \label{FormanRicciEdge}
{ \rm Ric}_{F,r}(e) = w_e \left( \frac{w_{v_1}}{w_e} +  \frac{w_{v_2}}{w_e}  - \sum_{e_{v_1}\ \sim\ e,\ e_{v_2}\ \sim\ e} \left[\frac{w_{v_1}}{\sqrt{w_e w_{e_{v_1} }}} + \frac{w_{v_2}}{\sqrt{w_e w_{e_{v_2} }}} \right] \right)\,.
\end{equation}

\begin{rem}
	One could make appeal to the reduced Forman curvature defined in Formula (\ref{FormanRicciEdge}) above as a ``proxy" for the computation of the ``true'' Forman curvature, especially in Imaging tasks, thus simplifying and speeding up computations. However, it is to be expected that the results obtained thus restricting oneself to the ``skeleton'' of the image will render only rather weak approximations of the results obtained when including, via Formula (\ref{eq:Forman-General}), the essential information contained in the 2-skeleton (pixels, voxels, etc.) of the given image (or mesh). 
	Preliminary experiments show that while the graph curvature approximates well the ``full" Forman one for natural images (see Figure 2 below), they are, as expected, more noisy for CT images. 
	Still, this approach appears to be most useful in settings where little geometric content is available, being replaced by probabilistic or even more general densities, typical cases being the analysis of MRI images or texture understanding.   
	
	\begin{figure}[htb]%\label{fig:HeightTorus}
		\centering
		\includegraphics[scale=0.55]{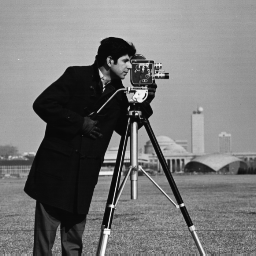}
		
		\includegraphics[scale=0.56]{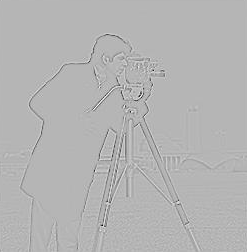} \includegraphics[scale=0.56]{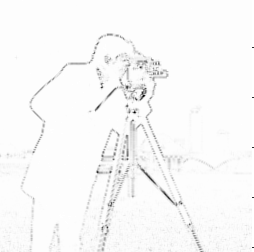}
		%{fr_sketch.jpg}
		\caption{{\bf Above:} A standard test image viewed as a complex of squares, the weight of each square (i.e. {\it pixel}), being equal to its grayscale level, and the weight of each edge being obtained via the Pythagorean Theorem from the grayscale levels, viewed as heights in the {\it stick model}. 
			{\bf Below left:}  The full Forman-Ricci curvature ${\rm Ric}_F$ of the considered complex. 
			{\bf Below right:} The reduced Forman-Ricci curvature ${\rm Ric}_{F,r}$ of the (geometric) dual graph, with vertex weights $w(u) = w(f)$, vertex $u$ being the center of the square $f$, and weight of dual edges equal those of the corresponding original ones (i.e. $w(e^*) = w(e), e^* \perp e$).}
	\end{figure}
\end{rem}

The  straightforward analogue of Formulas (2.) and (2.2) is then 
\begin{equation} \label{eq:RicciFlowNtwks-ch1}
\tilde{\gamma} (e)  - \gamma (e) = - {\rm Ric}_F (\gamma (e)) \cdot \gamma (e)\,;%
\end{equation}
where $\tilde{\gamma} (e)$ denotes the new (updated) value of $\gamma (e)$ after one time step.  Note that in our discrete setting, lengths are replaced by the (positive) 
edge weights. Also, for convenience and clarity, time is assumed to evolve in discrete steps, thus each ``clock" has a length of 1.  (However, in practice, this can and should be adapted to fit the specific circumstance of a given setup - see  also \cite{SSAZ}.)

\begin{rem}
 A derived scalar curvature and flow can also be considered (and easily derivable fro the Ricci curvature and flow), as shown in \cite{WJS},  Formulas (10) and (9), respectively.
\end{rem}

We should note that this simple and direct approach to the Ricci flow is excellently suited for applications where a short time flow is involved (see, e.g. \cite{ASZ11}, \cite{SSAZ}), it is not appropriate in situations where a long type flow should be employed, since it is not equipped with a fitting analogue for graphs of a surface of limit geometry, as in the smooth \cite{Ha1} and combinatorial \cite{CL} cases. 

However, for the case of higher dimensional complexes (i.e. for the original setting of Forman's definition), and in 
particular for polyhedral surfaces, such a limit geometry can be considered, due to Bloch's extension \cite{Bl} of Forman's
 work to allow for a fitting Gauss-Bonnet type theorem (which was an unfortunate lacuna in \cite{Fo} due to the algebraic approach of the original definition of Ricci curvature therein). This is possible due to his definition of Euler characteristic for {\it posets} (and, in particular for cell complexes), of which we do not bring the details here, but rather refer the reader to the original paper \cite{Bl}. % as well as to our forthcoming paper \cite{SS} on its implications on complex networks. 
 Suffice to say, that this is achieved by (somewhat technically) defining a triplet of combinatorial combinatorial curvatures, for each of the dimensions 0, 1 and 2, denoted by $R_1, R_2$ and $R_3$,  respectively. Since their specific definitions are quite technical, we do not detail them here, but content ourselves to noting that they satisfy a Gauss-Bonnet type theorem, more precisely one has

 \begin{thm}[Bloch, \cite{Bl}, Theorem 2.4] \label{thm:Bloch}
 	Let $X$ be 2-dimensional cell complex. Then:
 	
 	\begin{equation} \label{eq:Bloch-GaussBonnet1}
 	\sum_{v \in F_0}R_0(v) - \sum_{e \in F_1}R_1(e) + \sum_{f \in F_2}R_2(f) = \chi(X)\,.
 	\end{equation}
 \end{thm}
 
 The full relevance of the result above and of his importance for our study is the fact that, for polyhedral surfaces, the 1-dimensional curvature function and the Forman Ricci curvature coincide, that is 
 
  \begin{equation} \label{eq:ER1=Ric}
 R_1(e) = {\rm Ric}_F(e)\,.
 \end{equation} 
 
Given that polyhedral complexes represent a special case of cell complexes, from formulas (\ref{eq:Bloch-GaussBonnet1}) and  (\ref{eq:Bloch-GaussBonnet}) above imply that, for a 2-dimensional polyhedral complex, the following holds:

\begin{equation}  \label{eq:Bloch-GaussBonnet}
\sum_{v \in F_0}R_0(v) - \sum_{e \in F_1}{\rm Ric_F}(e) + \sum_{f \in F_2}R_2(f) = \chi(X)\,.
\end{equation}

Bloch's result above not only patches the serious aforementioned gap in Forman's work regarding his discretization of Ricci curvature, it also enables us to consider a proper long time Ricci flow, by allowing for the definition of  \textit{prototype} (or \textit{reference}, \textit{model}) \textit{networks}, as minimal (in the sense of having the minimal number of 2-faces) 2-dimensional polyhedral complexes with the Euler characteristics prescribed by the Gauss-Bonnet formula (\ref{eq:Bloch-GaussBonnet}). In consequence, one can define a network to be {\em spherical}, {\em Euclidean} or {\em hyperbolic}, if its Euler characteristic is ``$>$'', ``$=$" or ``$<$" %$\gtreqqless$
0, respectively. More precisely, we can thus propose the following 

\begin{defn}[Prototype networks]
	Let X be a 2-dimensional polyhedral complex with Euler characteristic $\chi$ as given by the Gauss-Bonnet formula. Then we define $\chi$ to be
	\begin{enumerate}
		\item Spherical, if $\chi >0$;
		\item Euclidean, if $\chi = 0$;
		\item Hyperbolic, if $\chi <0$.
	\end{enumerate}
\end{defn}

However, it is more convenient and efficient to consider, in analogy with the classical (surfaces) case, an alternative, equivalent definition of prototype complexes, as being complexes of constant curvature, i.e. such that $R_1 = {\rm const.}$, thus replacing conditions (1)-(3) above with the following ones, respectively:

(1') $R_1 > 0$;

(2')  $R_1 = 0$;

(3')  $R_1 < 0$.

This alternative definition is not only more simple, it also facilitates the definition of the desired long term flow, which we formulate here only for the case of 2-dimensional polyhedral complexes. More precisely, we define -- in analogy with the classical case \cite{Ha1} (see also (Equation (2.3) above)), and keeping in mind that in this case $R_1(e) = {\rm Ric}_F(e)$ -- the {\it Normalized Forman-Ricci flow} for networks  2-dimensional polyhedral complexes), by

\begin{eqnarray} \label{eq:NormalizedFormanRicciFlow}
\frac{\partial \gamma(e)}{\partial t} = -\left( \rm{Ric_F} \left( \gamma(e) \right)  - \rm{\overline{Ric}_F}\right) \cdot \gamma(e) \; ,
\end{eqnarray}
where $\rm{\overline{Ric}_F}$ denotes the {\it mean} Forman-Ricci curvature.

It follows therefore that, by using Bloch's results, we are able to define a proper notion of background (limit) geometry for complexes that allows for the study of long term evolution of 2-dimensional polyhedral complexes. In consequence, it enables us to study their (topological) complexity, dispersion of geodesics, volume growth, etc.; that is all the essential properties captured by its geometric type. 
To be sure, the next step would be to experimentally explore  how the above defined flow compares to similar geometric flows, starting with the best known one, namely the combinatorial flow. The experiments conducted so far suggest that, indeed, the notion of limit geometry introduced here is a viable one. 
However, we have still both to provide theoretical proofs of the basic convergence theorems similar to those in \cite{Ha1}, \cite{CL}, as well as to perform large scale experiments, these being tasks that for now we defer for future study.

%---------------------------------------------------

%------------------------------------------------------------

\section*{Acknowledgement}
The author would like to thank his former and present students who helped produce the figures in the text, to J{\"u}rgen Jost for his warm hosting at the Max Planck Institute for Mathematics in the Sciences, Leipzig, where a substantial part of this paper was conceived and written, to David Xianfeng Gu, whose  keen %enthusiastic
interest and warm encouragement started this project. Thanks are also due to Feng Luo and Areejit Samal for useful discussion regarding the Forman-Ricci flow.

%The author would like to thank David Xianfeng Gu for his keen interest, his stimulating questions and for %his patience in
%guiding him through the implementation aspects of the combinatorial Ricci flow and to Feng Luo for his constructive critique and his encouragement.

%Special thanks are due to the anonymous reviewer of an earlier version of this paper, for his guiding, stimulating questions that corrected and clarified the authors' own thoughts, and for the new, important bibliographical references that focused the work on better directions (both for the present as well as for further refinements of this study).
%
%Finally, the author would like to thank his former students who helped produce some of the figures in the text, as well as to Eli Appleboim for his assistance in this matter.

%%%%%%%%%%%%%%%%%%%%%%%%%%%%%%%%%%%%%%%%%%%%%%

%%%%%%%%%%%%%%%%%%%%%%%%%%%%%%%%%

\end{document}